\documentclass{aptpub}

\usepackage{latexsym} 
\usepackage{amssymb}
\usepackage{cite}
\usepackage{amsmath}
\usepackage{hyperref}
\usepackage{bbm}
\usepackage{mathtools}
\newtheorem{algorithm}[theorem]{Algorithm}

\authornames{David Aristoff} 
\shorttitle{An ergodic theorem for weighted ensemble} 

\begin{document}

\title{An Ergodic Theorem for the Weighted Ensemble Method} 

\authorone[Colorado State University]{David Aristoff} 

\addressone{841 Oval Drive, 
Fort Collins, CO 80523} 
\emailone{aristoff@math.colostate.edu} 

\begin{abstract}
We study weighted ensemble,  
an interacting particle method for 
sampling distributions of Markov chains 
that has been used in computational 
chemistry since the 1990s. Many important applications of weighted ensemble require the computation of long time averages. We establish the consistency of 
weighted ensemble in this setting by proving an ergodic 
theorem for time averages. 
As part of the proof, we derive explicit variance formulas that could be useful for optimizing the method.
\end{abstract}

\keywords{Interacting particle system, sequential Monte Carlo, 
ergodic theorem, weighted ensemble, importance sampling, steady state, stratification}

\ams{65C05}{65C35;65C40}

         \section{Introduction}

{Weighted ensemble}~\cite{
bhatt2010steady,
chong2017path,costaouec2013analysis,
darve2013computing,dickson,donovan2013efficient,
huber1996weighted,
rojnuckarin1998brownian,
rojnuckarin2000bimolecular,zhang2007efficient,
zhang2010weighted,
zwier2015westpa} is an importance sampling
method, based on interacting particles, for distributions associated 
with a Markov chain. 
In this article, we 
focus on sampling the average 
of a function
with respect to the steady state distribution of 
a generic Markov chain. By generic, 
we mean that the only thing 
we might know about the Markov chain is 
how to sample it; in particular, 
we may not know its stationary 
distribution up to a 
normalization factor.

Weighted 
ensemble consists of a collection 
of evolving {\em particles} with 
associated {\em weights}. In this sense, weighted ensemble can be understood 
as a kind of sequential Monte 
Carlo method~\cite{doucet2001sequential,del2004feynman,del2014particle,del2005genealogical,weare,assaraf,doucetSMC, webber2,webber}.
In weighted ensemble, 
the particles evolve between 
{\em selection} steps according 
to the law of the underlying Markov 
chain. In each selection step, 
some of the particles are {copied} 
while others are {killed}; 
the resulting particles are 
given new weights so that  weighted 
ensemble is statistically 
{\em unbiased}~\cite{zhang2010weighted}.

The selection step is based on 
dividing the particles into 
{\em bins}, where  
the particles in each bin are 
resampled according to their 
relative weights. 
In practice, the binning, 
and the number of copies 
in each bin, 
should be chosen so that important 
particles survive and irrelevant particles 
are killed. 
The definition 
of the bins, and how many copies to 
maintain in each, 
requires some care.  
With appropriate choices, weighted 
ensemble can have drastically smaller 
variance than direct Monte 
Carlo, or independent particles; see 
the references above, or~\cite{zuckermanchong} for a 
more complete list.

Weighted ensemble was developed 
for applications in computational 
chemistry~\cite{huber1996weighted} 
ranging from state space 
exploration~\cite{dickson} to protein association~\cite{huber1996weighted} and protein 
folding~\cite{zwier2010reaching}. One important  
application we have in mind is 
the computation of the mean
time for a protein to unfold~\cite{bhatt2010steady}. 
This time can be reformulated as the inverse of the steady state 
flux into the unfolded state of the underlying Markovian 
molecular dyamics, with an added sink 
in the unfolded state and source in the 
folded state~\cite{hill}. This dynamics can approach its steady state on time scales significantly smaller than the mean unfolding time~\cite{zuckermanweb}.
As the flux into the unfolded state is usually very small,  
importance sampling is 
needed to 
estimate it 
with substantial precision~\cite{bhatt2010steady,suarez}.

Other unbiased methods, differing 
from weighted ensemble in that they 
usually sample finite-time quantities rather than 
ergodic averages, 
include Adaptive Multilevel 
Splitting~\cite{brehier_AMS,brehier_AMS2,tony_AMS}, Forward Flux Sampling~\cite{allen2006forward}, 
and some sequential Monte Carlo methods~\cite{del2005genealogical,chraibi2018optimal,webber}.
This unbiased property allows 
for a relatively straightforward study of variance using 
martingale techniques~\cite{aristoff2016analysis,
brehier_AMS2,del2004feynman,del2005genealogical}. 
In this article, we extend these techniques to study the long-time stability of weighted 
ensemble.

Our main contribution here is a proof of
the consistency of
weighted ensemble 
via an ergodic theorem. We believe that this is the first ergodic theorem 
for an interacting particle system in which 
the interactions come from resampling.

A secondary contribution 
comes from explicit formulas for the 
variance of weighted ensemble at finite particle number. The proof 
of the ergodic theorem is a 
straighforward consequence of these formulas.
On the theoretical side, our 
variance formulas are handy for understanding 
the rate of weighted ensemble convergence, and on the 
practical side, they
could be used for optimizing the 
method. We mostly leave this discussion 
to other works, including
our companion paper~\cite{aristoff2018steady}; 
see also~\cite{aristoff2016analysis}
and the references above.

This article is organized as follows. In 
Section~\ref{sec:description}, we describe 
weighted ensemble in detail. In Section~\ref{sec:ergodic_theorem}, we state our 
main results, including the unbiased property (Theorem~\ref{thm:unbiased}), the ergodic theorem (Theorem~\ref{thm_main}), and the variance formulas (Theorem~\ref{thm:variance_expressions}). In Section~\ref{sec:compare}, we compare weighted 
ensemble to direct Monte Carlo, and give a simple example illustrating the potential gain. 
All of our proofs are in Section~\ref{sec:derivations}.

\section{Description of the method}\label{sec:description}

Weighted ensemble consists of 
a fixed number, $N$, of 
particles 
belonging to 
a common state space, each 
carrying a positive scalar weight, and undergoing 
repeated selection and mutation steps. 
In the selection step, some of the particles are copied, and others are killed, according to a 
stratification or {\em binning} scheme. In the mutation step, the particles 
evolve according to an underlying Markov kernel $K$.

At time $t$ before selection, the particles, called {\em parents}, are $\xi_t^1,\ldots,\xi_t^N$.
At time $t$ after selection, the particles, called {\em children}, are  
$\hat\xi_t^1,\ldots,\hat\xi_t^N$. The 
weights of the parents and children are 
$ \omega_t^1,\ldots,\omega_t^N$ and
$\hat\omega_t^1,\ldots,\hat\omega_t^N$, respectively.
The following diagram illustrates weighted ensemble evolution:
\begin{align}\begin{split}\label{eq:sel_and_mut}
\begin{array}{ccc}
\text{parents }\{\xi_t^i\}^{i=1,\ldots,N} & \text{parents' weights } \{\omega_t^i\}^{i=1,\ldots,N} & \textup{user-chosen}\\
\downarrow \text{selection} & \downarrow \text{selection} & \textup{parameters} \\
\text{children }\{\hat{\xi}_t^i\}^{i=1,\ldots,N} & \text{children's weights } \{\hat{\omega}_t^i\}^{i=1,\ldots,N} & \leftarrow\\
\downarrow \text{mutation}& \downarrow \text{mutation}\\ 
\text{new parents }\{\xi_{t+1}^i\}^{i=1,\ldots,N} & \text{new parents' weights } \{\omega_{t+1}^i\}^{i=1,\ldots,N}
\end{array}\end{split}
\end{align}

The initial particles 
$\xi_0^1, \ldots,\xi_0^N$ 
can be arbitrary. The initial 
weights must be strictly 
positive and sum to one: $\omega_0^i>0$ for all $i$, 
and $\omega_0^1+\ldots + \omega_0^N =1$. The children are initially just copies of their parents, but they evolve forward in time conditionally independently. When we say {\em conditionally independent}, 
we mean conditional on the $\sigma$-algebra 
representing the information 
from~\eqref{eq:sel_and_mut} up to the current time.

Weighted ensemble requires the 
user to choose, before the selection 
step at time $t$, a 
collection of
nonempty
{\em bins} that partition 
the set of parents, as well as a  
{\em particle allocation} that 
defines the number of children in each bin. 
We write $u$ for the bins, 
and $N_t(u) \ge 1$ for the number of children 
in bin $u$ at time $t$. We require that $\sum_u N_t(u) = N$. The bins 
can change in time, but for simpler 
notation we leave this implicit.  

In the selection step, the children in each 
bin
are obtained by sampling with replacement from the 
parents in the bin, according to their 
weight distribution, as many times as the particle allocation specifies. The children's weights in 
each bin
are all the same after selection, and the 
total weight in the bin is preserved~\cite{costaouec2013analysis}. 

In more detail, define the total weight in bin $u$ at time $t$ as 
\begin{equation}\label{eq:bin_weights}
\omega_t(u) = \sum_{i:\xi_t^i \in u} \omega_t^i.
\end{equation}
The numbers, $N_t^i$, of children of the parents $\xi_t^i \in u$ are conditionally
multinomial: 
\begin{equation}\label{eq:multinomial}
\{N_t^i:\xi_t^i \in u\} \sim \textup{Multinomial}\left(N_t(u), \left\{\frac{\omega_t^i}{\omega_t(u)}:\xi_t^i \in u\right\}\right).
\end{equation}
Children are assigned to the same bins as their parents, with weights
\begin{equation}\label{eq:hat_weights}
\hat{\omega}_t^j = \frac{\omega_t(u)}{N_t(u)},\qquad \text{if }\hat{\xi}_t^j  \in u.
\end{equation}
Selections in distinct 
bins are conditionally independent.

In the mutation step, 
the children evolve 
conditionally independently
via $K$:
\begin{equation}\label{eq:mutation}
(\xi_{t+1}^1,\ldots,\xi_{t+1}^N) \sim K(\hat{\xi}_t^1,\cdot) \times \ldots \times K(\hat{\xi}_t^N,\cdot).
\end{equation}
The 
weights do not change during the 
mutation step. Thus
\begin{equation}\label{eq:weights}
{\omega}_{t+1}^{j} = \hat{\omega}_{t}^j, \qquad j=1,\ldots,N.
\end{equation}

We summarize weighted ensemble in the 
following algorithm. 
\begin{algorithm}\label{alg:WE}
Choose initial weights $\omega_0^1,\ldots,\omega_0^N>0$ summing to $1$ and initial particles $\xi_0^1,\ldots,\xi_0^N$. Then iterate over $t \ge 0$:

\begin{itemize}
\item[]{\em Selection step}

\item[1.] Partition the parents  $\xi_t^1,\ldots,\xi_t^N$ into a collection of {bins}. 
\item[2.] Assign a number $N_t(u)\ge 1$ of children to the parents in each bin $u$. 
\item[3.] Sample $N_t(u)$ children from the parents in bin $u$, with replacement, using
$$Pr(\textup{sample }\xi_t^j \text{ in bin }u) = \frac{\omega_t^j}{\omega_t(u)},$$
where $\omega_t(u) := \sum_{i\,:\,\xi_t^i \in u}\omega_t^i$ is the total weight in bin $u$.
\item[4.]
Give all the children 
in bin $u$ the same weight
$$\hat{\omega}_t^j = \frac{\omega_t(u)}{N_t(u)}, \qquad \text{if }\hat{\xi}_t^j \in u.$$
\end{itemize}

\begin{itemize}
\item[]{\em Mutation step}

\item[5.] Evolve the 
children $\hat{\xi}_t^1,\ldots,\hat{\xi}_t^N$ conditionally independently 
using $K$ to get the 
next parents $\xi_{t+1}^1,\ldots,\xi_{t+1}^N$.
Keep the same weights 
$\omega_{t+1}^j = \hat{\omega}_t^j$, $j=1,\ldots,N$.

\end{itemize}

\end{algorithm}

\subsection{Algorithm details and remarks}

A few remarks are in order to clarify Algorithm~\ref{alg:WE}.
\begin{itemize}
\item
We abuse notation by writing $\xi_t^i \in u$ or $\hat{\xi}_t^j \in u$ to indicate that $\xi_t^i$ or $\hat{\xi}_t^j$ is in bin $u$, even though the bins form a partition of the particles, and not (necessarily) a partition of the state space.
\item
A child is simply a copy of its parent: if $\hat{\xi}_t^j$ is a child of $\xi_t^i$, then $\hat{\xi}_t^j = \xi_t^i$. The indices of the children are not important, 
so specifying the number of children of each parent is enough to define them.
\item
Since the weights don't change in the mutation step and the the selection step preserves the total weight, the total weight is constant in time: $\omega_t^1 + \ldots +  \omega_t^N = 1$ for all $t \ge 0$. We 
discuss the importance of this 
in Remark~\ref{rmk:total_weight}.
\item
We assume that the bins and particle allocation at time $t$ are included in  ${\mathcal F}_t$,
the $\sigma$-algebra generated 
by the information from Algorithm~\ref{alg:WE} just before the 
$t$-th selection step. We also write $\hat{\mathcal F}_t$ 
for the $\sigma$-algebra generated 
by the information from Algorithm~\ref{alg:WE} just after the 
$t$-th selection step. See Section~\ref{sec:notation} for details.
\item We assume multinomial sampling in the bins because it leads to
simple explicit variance expressions in terms 
of intrabin variances. In Remark~\ref{rmk:residual}, we comment on residual 
sampling, which performs much better than 
multinomial resampling and still admits nice
variance formulas.
\end{itemize}

For our ergodic theorem, the bins and particle 
allocation can be arbitrary. To
actually do better than
direct Monte Carlo, they must be judiciously 
chosen. The most common strategy is to define 
bins based on a carefully constructed partition 
of state space -- particles occupy the 
same bin when they belong to the same 
element of the partition -- and then 
allocate children approximately uniformly among 
these bins. Some knowledge about the underlying 
problem is needed to choose the bins, 
but this strategy has had considerable success,  
as the references in the introduction attest (see~\cite{zuckermanchong} for a 
mostly current list of application papers).
We propose a different strategy in our companion paper~\cite{aristoff2018steady} that uses 
our variance analysis below. We summarize that strategy in Section~\ref{sec:compare} below.

\section{Main results}\label{sec:ergodic_theorem}

\subsection{Unbiased property}

We begin with the unbiased property of weighted ensemble. This property was previously noted in~\cite{zhang2010weighted},
and proved in a slightly different setting in~\cite{aristoff2016analysis}.

\begin{theorem}[Unbiased property]\label{thm:unbiased}
For each $t\ge 0$ and all bounded measurable $g$, 
\begin{equation*}
{\mathbb E}\left[\sum_{i=1}^N \omega_t^i g(\xi_t^i)\right] = \int K^t g\,d\nu,
\end{equation*}
where $\nu$ is the weighted ensemble initial distribution, $\int g\,d\nu := {\mathbb E}[\sum_{i=1}^N \omega_0^i g(\xi_0^i)]$.
\end{theorem}

We could interpret Theorem~\ref{thm:unbiased} as follows. If 
$(X_t)$ is a Markov chain with 
kernel $K$ and initial 
distribution $\nu$, then 
${\mathbb E}\left[\sum_{i=1}^N \omega_t^i g(\xi_t^i)\right] = {\mathbb E}[g(X_t)]$. 
In this sense, weighted ensemble 
gives unbiased estimates of 
the law of the underlying Markov chain.

\subsection{Ergodic theorem}

To ensure that weighted ensemble is ergodic, the underlying 
Markov kernel $K$ must be ergodic in some sense. We assume that $K$ is uniformly ergodic~\cite{douc}:

\begin{assumption}\label{A1}
There is $c>0$, $\lambda \in [0,1)$ and a probability measure $\mu$ such that $$\|K^t(\xi,\cdot) - \mu(\cdot)\|_{TV} \le c \lambda^t, \qquad \text{for all }\xi\text{ and all }t\ge 0.$$
\end{assumption}

Here and below, $f$ is a fixed bounded measurable function. 
 
\begin{theorem}[Ergodic theorem]\label{thm_main}
If Assumption~\ref{A1} holds, then with probability $1$,
\begin{equation}\label{eq:ergodic}
\lim_{T \to\infty} \frac{1}{T}\sum_{t=0}^{T-1}\sum_{i=1}^N \omega_t^i f(\xi_t^i)  =\int f\,d\mu.
\end{equation}
\end{theorem}

Convergence of the {\em mean} of the time average in~\eqref{eq:ergodic}, at the same rate as direct Monte Carlo, follows from Assumption~\ref{A1} and the unbiased property (Theorem~\ref{thm:unbiased}). For the ergodic theorem to hold, 
and for weighted ensemble to beat 
direct Monte Carlo, the {\em variance} of the time average should be sufficiently small. Well-behaved variance is not automatic for unbiased methods; see Remark~\ref{rmk:total_weight} below.

\subsection{Variance formulas}

Here, we give exact, finite $N$ formulas for the variance of weighted ensemble, based on a martingale 
decomposition. To get nice concise formulas, 
we need some notation. Define the intrabin distributions 
\begin{equation*}
\eta_t^u = \sum_{i:\xi_t^i \in u} \frac{\omega_t^i}{\omega_t(u)} \delta_{\xi_t^i},
\end{equation*}
where $\delta_{\xi}$ is the Dirac delta distribution centered at $\xi$. 
Define also
\begin{equation*}
h_{t,T} = \sum_{s=0}^{T-t-1} K^sf.
\end{equation*}
For a probability measure $\eta$ and bounded measurable function $g$, define 
\begin{equation*}
\eta(g) = \int g\,d\eta, \qquad\textup{Var}_\eta(g) = \eta(g^2)-\eta(g)^2,
\end{equation*}
and in particular, let $\textup{Var}_K g(\xi) := \textup{Var}_{K(\xi,\cdot)}(g) = Kg^2(\xi) - (Kg)^2(\xi)$.

\begin{theorem}[Variance formulas]\label{thm:variance_expressions}
For each time $T >0$,  
\begin{align}
&\textup{Var}\left(\frac{1}{T}\sum_{t=0}^{T-1}\sum_{i=1}^N \omega_t^if(\xi_t^i) \right) 
\label{eq:WE_var}\\
&= \frac{1}{T^2}\textup{Var}\left(\sum_{i=1}^N \omega_0^i h_{0,T}(\xi_0^i)\right) \qquad \textup{(initial condition variance)} \label{eq:var_ic}\\
& + \frac{1}{T^2}\sum_{t=0}^{T-2}{\mathbb E}\left[\sum_u \frac{\omega_t(u)^2}{N_t(u)}\textup{Var}_{\eta_t^u}(Kh_{t+1,T})\right]
\qquad \textup{(selection variance)}
\label{eq:var_sel} \\
& \,\,\,+ \frac{1}{T^2}\sum_{t=0}^{T-2} {\mathbb E}\left[\sum_u\frac{\omega_t(u)^2}{N_t(u)}\eta_t^u(\textup{Var}_K h_{t+1,T})\right] \qquad \textup{(mutation variance)}.\label{eq:var_mut}
\end{align}
\end{theorem}

The expression in~\eqref{eq:var_ic} can be interpreted as the 
variance coming from the initial condition, while the expressions in~\eqref{eq:var_sel} and~\eqref{eq:var_mut} can be understood as the 
variances arising from each selection 
and mutation step, respectively. 

Using Theorem~\ref{thm:variance_expressions}, the proof of the ergodic theorem is straightforward. 
Under Assumption~\ref{A1}, we can show that variances of $h_{t,T}$ and $Kh_{t,T}$ are uniformly bounded in $t$ and $T$. This makes the weighted ensemble variance $\textup{Var}\left(\frac{1}{T}\sum_{t=0}^{T-1}\sum_{i=1}^N \omega_t^if(\xi_t^i) \right) = O(1/T)$. 
The ergodic theorem then
follows from standard arguments. 
Note that the variance expressions~\eqref{eq:WE_var}-\eqref{eq:var_mut} by themselves 
do not require Assumption~\ref{A1}.

Beyond the ergodic theorem, these variance formulas are interesting 
in their own right, since they could be 
used to design binning and particle allocation schemes that minimize the weighted ensemble variance. 
Indeed, that is what we have done in our companion paper~\cite{aristoff2018steady} (see also~\cite{aristoff2016analysis}). We discuss this more in the next section.

\section{Comparison to direct Monte Carlo}\label{sec:compare}

There are many existing works showing that 
weighted ensemble can provide significant gains over direct Monte Carlo: see for instance the references list in~\cite{zuckermanchong}, 
and our companion paper~\cite{aristoff2018steady}. 
The main goal of this article is to prove the consistency of weighted ensemble from explicit variance formulas, and
not to reaffirm this point. We include, however, a brief 
discussion here.  

Weighted ensemble 
works by reducing the mutation variance 
\eqref{eq:var_mut}, compared to that of 
direct Monte Carlo (see~\eqref{eq:DMC_mutvar} below), 
via the selection step. This comes however at the cost of a positive selection variance 
\eqref{eq:var_sel}, compared to direct Monte Carlo which has selection variance 
equal to zero. Thus, 
weighted ensemble beats direct Monte Carlo if  
the reduction in mutation variance is greater than 
the selection variance cost. (The 
initial condition variances may
be ignored, since they are $O(1/T^2)$ for weighted 
ensemble and direct Monte Carlo, while 
the overall variances are $O(1/T)$.)

In more detail, the weighted ensemble mutation variance is
\begin{equation}\label{eq:WE_mutvar}
\frac{1}{T^2}\sum_{t=0}^{T-2}{\mathbb E}\left[\sum_u\frac{\omega_t(u)^2}{N_t(u)}\eta_t^u(\textup{Var}_Kh_{t+1,T})\right].
\end{equation}
A Lagrange multiplier calculation shows that the 
expression inside the expectation in~\eqref{eq:WE_mutvar} is minimized over $N_t(u)$, subject to the constraint $\sum_u N_t(u) = N$, when 
\begin{equation}\label{eq:Ntu}
N_t(u) \approx \frac{N\omega_t(u)\sqrt{\eta_t^u(\textup{Var}_Kh_{t+1,T})}}{\sum_u \omega_t(u)\sqrt{\eta_t^u(\textup{Var}_Kh_{t+1,T})}}.
\end{equation} 
Plugging this into~\eqref{eq:WE_mutvar} 
makes the weighted ensemble mutation variance 
\begin{equation}\label{eq:WE_mutvar2}
\approx \frac{1}{NT^2}\sum_{t=0}^{T-2}{\mathbb E}\left[ \left(\sum_u\omega_t(u)\sqrt{\eta_t^u(\textup{Var}_Kh_{t+1,T})}\right)^2\right].
\end{equation}
By Jensen's inequality,~\eqref{eq:WE_mutvar2} is less than or equal to 
\begin{equation}\label{eq:direct_MC_variance}
\frac{1}{NT^2} \sum_{t=0}^{T-2}{\mathbb E}\left[\sum_u \omega_t(u)\eta_t^u(\textup{Var}_Kh_{t+1,T})\right], 
\end{equation}
which, by the unbiased property of weighted ensemble (Theorem~\ref{thm:unbiased}), is 
exactly the direct Monte Carlo mutation 
variance (see Remark~\ref{rmk:DMC} below). 
This 
form of optimal mutation variance gain was originally 
observed in~\cite[Remark 4.1]{aristoff2016analysis}.

The selection variance~\eqref{eq:var_sel} is small whenever, at each time $t$, bins 
$u$ are chosen inside which $Kh_{t+1,T}$ does 
not vary too much. With enough particles 
and bins, it is possible to keep the selection variance arbitrarily small, while 
also controlling the mutation variance 
by keeping the particle allocation close to 
the optimal~\eqref{eq:Ntu}. Even with modest numbers 
of particles and bins, weighted ensemble has 
proven useful; see the list of applications in~\cite{zuckermanchong}, most of which 
use relatively small $N$.

In our companion 
paper~\cite{aristoff2018steady}, we propose an optimization strategy based on choosing 
the particle allocation to minimize mutation variance, and the bins to minimize selection variance. There, the particle allocation is a simplified version of~\eqref{eq:Ntu}, in which we 
use the limit $h:= \lim_{T\to \infty} (h_{t,T}- (T-t)\int f\,d\mu)$ in place of $h_{t+1,T}$, 
and the bins are chosen to make the intrabin 
variances of $Kh$ small. 
Of course, estimating $Kh$ and $\textup{Var}_Kh$ is a difficult 
problem. In~\cite{aristoff2018steady} we propose 
using a Markov state model 
to get ``cheap'' approximations of $Kh$ and $\textup{Var}_Kh$; such models 
are already commonly used for 
preconditioning weighted 
ensemble simulations~\cite{bhatt2010steady,copperman1,copperman2}.

\subsection{Example}

Below is a simple example 
illustrating the variance reduction in weighted ensemble, compared to direct Monte Carlo, in the context 
of our variance analysis above.
Consider a Markov chain on $3$ states, with transition matrix
$$K(1,2) =  K(2,3) = \delta,\qquad  
K(1,1) = K(2,1) = 1-\delta, \qquad
K(3,1) = 1,$$ where 
$\delta>0$ is small, and we take $f(1) = f(2) = 0$  and $f(3) = 1$.

For weighted ensemble, we will assign 
particles to the same bin if and only if they occupy the same point in space. Following 
the usual method in the literature, we allocate an approximately 
equal number of children to each bin.

\begin{figure}
\includegraphics[width=13cm]
{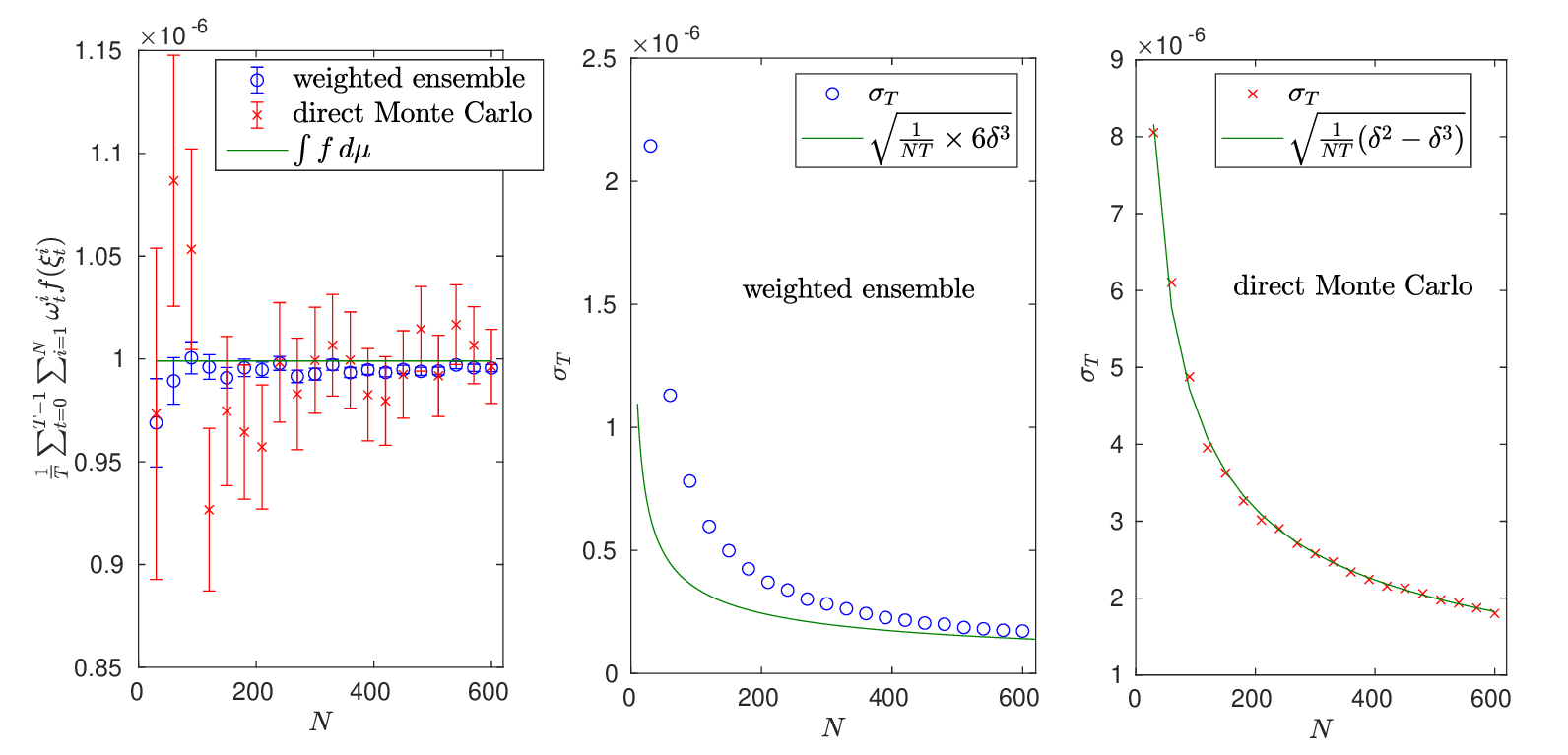}
\caption{Comparison of weighted 
ensemble with direct Monte Carlo when $\delta = 0.001$ and $T = 500$. {\em Left}: Average values of $\frac{1}{T}\sum_{t=0}^{T-1}\sum_{i=1}^N \omega_t^i f(\xi_t^i)$ vs. $N$, computed from $10^4$ independent trials, for weighted ensemble and direct Monte Carlo. Error 
bars widths are $\sigma_T/10^2$ where $\sigma_T^2$ are 
the empirical variances. {\em Center}: 
Weighted ensemble empirical standard deviation compared with~\eqref{eq:var_WE_example}. {\em Right}: 
Direct 
Monte Carlo empirical standard deviation compared with~\eqref{eq:var_DMC_example}.}
\label{fig1}
\end{figure}

In this case, the variance of weighted ensemble can be estimated as
\begin{equation}\label{eq:var_WE_example}
 \frac{1}{T^2}\sum_{t=0}^{T-2}{\mathbb E}\left[\sum_{i=1}^3\frac{\omega_t(i)^2}{N_t(i)}\textup{Var}_K h_{t+1,T}(i)\right] 
\approx \frac{1}{T}\sum_{i=1}^3 \frac{\mu(i)^2}{N/3}\textup{Var}_K h_{t+1,T}(i) 
\approx \frac{6\delta^3}{NT},
\end{equation}
where the first approximation in~\eqref{eq:var_WE_example} holds for large enough $N$ and $T$, and the second approximation uses direct calculations, dropping terms of higher order than $\delta^3$. The variance of direct Monte Carlo (see Remark~\ref{rmk:DMC}) can be estimated by 
\begin{equation}\label{eq:var_DMC_example}
 \frac{1}{NT^2}\sum_{t=0}^{T-2}\nu(K^t(\textup{Var}_Kh_{t+1,T})) 
\approx \frac{1}{NT}\mu(\textup{Var}_K h_{t+1,T})
\approx \frac{\delta^2-\delta^3}{NT}
\end{equation}
for large $T$. Figure~\ref{fig1} shows numerical confirmation of these estimates.
Note the smaller variance (higher order in $\delta$) for weighted ensemble, compared to direct Monte Carlo.

\section{Derivations}\label{sec:derivations}

\subsection{Notation}\label{sec:notation}

Below, ${\mathcal F}_t$ is the $\sigma$-algebra generated 
by the parents and their weights from times $0 \le s \le t$, the 
children and their weights from 
times $0 \le s \le t-1$, and the bins and particle allocations from times $0 \le s \le t$. 
Meanwhile, $\hat{\mathcal F}_t$ is the 
$\sigma$-algebra generated by ${\mathcal F}_t$ together with the children and their weights at time $t$. Throughout, $g$ 
denotes a bounded measurable function, 
and $c$ a positive constant whose value 
can change between different equations. We will use the notation
\begin{equation}\label{eq:parents}
\textup{par}(\hat{\xi}_t^j) = \xi_t^i \Longleftrightarrow \xi_t^i \text{ is the parent of }\hat{\xi}_t^j.
\end{equation}

\subsection{One-step means}

\begin{lemma}\label{lem:one_step}
For each $i=1,\ldots,N$ and $t \ge 0$,
\begin{equation}\label{eq:mut_mean}
{\mathbb E}\left[\left.\omega_{t+1}^i g(\xi_{t+1}^i)\right|\hat{\mathcal F}_t\right] =  \hat{\omega}_t^i Kg(\hat{\xi}_t^i).
\end{equation}
At each time $t \ge 0$, for each bin $u$,
\begin{equation}\label{eq:bin_mean}
{\mathbb E}\left[\left.\sum_{i:\hat{\xi}_t^i \in u}\hat{\omega}_t^i g(\hat{\xi}_t^i)\right|{\mathcal F}_t\right] = 
\sum_{i:\xi_t^i \in u} \omega_t^i g(\xi_t^i).
\end{equation}
\end{lemma}

\begin{proof}
By~\eqref{eq:mutation}, ${\mathbb E}[g(\xi_{t+1}^i)|\hat{\mathcal F}_t] = Kg(\hat{\xi}_t^i)$. Now~\eqref{eq:mut_mean} follows from this and~\eqref{eq:weights}. 
Meanwhile, by~\eqref{eq:multinomial}, ${\mathbb E}[N_t^i|{\mathcal F}_t] = N_t(u)\omega_t^i/\omega_t(u)$ if $\xi_t^i \in u$. Thus by~\eqref{eq:hat_weights} and~\eqref{eq:parents},
\begin{align*}
{\mathbb E}\left[\left.\sum_{i:\hat{\xi}_t^i \in u}\hat{\omega}_t^i g(\hat{\xi}_t^i)\right|{\mathcal F}_t\right] 
&= \sum_{i:\xi_t^i \in u}{\mathbb E}\left[\left.\sum_{j:\textup{par}(\hat{\xi}_t^j) = \xi_t^i}
\hat{\omega}_t^j g(\hat{\xi}_t^j)\right|{\mathcal F}_t\right] \\
&= \sum_{i:\xi_t^i \in u}\frac{\omega_t(u)}{N_t(u)}g(\xi_t^i){\mathbb E}[N_t^i |{\mathcal F}_t] = \sum_{i:\xi_t^i \in u}\omega_t^i g(\xi_t^i).
\end{align*}
\end{proof}

\begin{lemma}[One-step means]\label{lem:one_step2}
For each time $t \ge 0$,
\begin{align}
&{\mathbb E}\left[\left.\sum_{i=1}^N \omega_{t+1}^i g(\xi_{t+1}^i)\right|\hat{\mathcal F}_t\right] = \sum_{i=1}^N \hat{\omega}_t^i Kg(\hat{\xi}_t^i), \label{eq:one_step1} \\
&{\mathbb E}\left[\left.\sum_{i=1}^N \hat{\omega}_t^i g(\hat{\xi}_t^i)\right|{\mathcal F}_t\right] = \sum_{i=1}^N  {\omega}_t^i g({\xi}_t^i). \label{eq:one_step2}
\end{align}
\end{lemma}
\begin{proof} This follows immediately from Lemma~\ref{lem:one_step}, by summing 
over the particles in~\eqref{eq:mut_mean} to get~\eqref{eq:one_step1}, and 
summing over the bins in~\eqref{eq:bin_mean} to get~\eqref{eq:one_step2}.
\end{proof}

\subsection{Proof of the unbiased property}

\begin{proof}[Proof of Theorem~\ref{thm:unbiased}]
By repeated application of Lemma~\ref{lem:one_step2} with the tower property,
\begin{equation}\label{eq:rep_tower}
{\mathbb E}\left[\left.\sum_{i=1}^N \omega_{t}^i g(\xi_{t}^i)\right|{\mathcal F}_0\right] = \sum_{i=1}^N \omega_0^i K^t g(\xi_0^i).
\end{equation}
Taking expectations in~\eqref{eq:rep_tower} gives the result.
\end{proof}

\subsection{Doob martingale and variance decomposition}

Below, define the Doob martingale
\begin{equation*}
D_t = {\mathbb E}\left[\left.\sum_{s=0}^{T-1}\sum_{i=1}^N \omega_s^i f(\xi_s^i)\right|{\mathcal F}_t\right], \qquad \hat{D}_t = {\mathbb E}\left[\left.\sum_{s=0}^{T-1}\sum_{i=1}^N \omega_s^i f(\xi_s^i)\right|\hat{\mathcal F}_t\right].
\end{equation*}
\begin{proposition}\label{prop:doob}
For $0 \le t \le T-1$, we have
\begin{align}\begin{split}\label{eq:doob_martingale}
D_t &= \sum_{s=0}^{t} \sum_{i=1}^N \omega_s^i f(\xi_s^i) + \sum_{i=1}^N \omega_{t}^i Kh_{t+1,T}(\xi_{t}^i), \\
\hat{D}_t &= \sum_{s=0}^{t} \sum_{i=1}^N \omega_s^i f(\xi_s^i) + \sum_{i=1}^N \hat{\omega}_{t}^i Kh_{t+1,T}(\hat{\xi}_t^i).
\end{split}
\end{align}
\end{proposition}
\begin{proof}
This comes from Lemma~\ref{lem:one_step2} by repeated application of the tower property.
\end{proof}

\begin{proposition}[Martingale variance decomposition]\label{lem:doob}
For each $T>0$,
\begin{align}\begin{split}
&\textup{Var}\left(\frac{1}{T}\sum_{t=0}^{T-1}\sum_{i=1}^N \omega_t^i f(\xi_t^i)\right) \\
&= \underbrace{\frac{1}{T^2}\textup{Var}(D_0)}_{\textup{initial condition variance}}  +
\underbrace{\frac{1}{T^2}\sum_{t=0}^{T-2}{\mathbb E}\left[(\hat{D}_t-D_t)^2\right]}_{\textup{selection variance}} +  \underbrace{\frac{1}{T^2}\sum_{t=0}^{T-2}{\mathbb E}\left[(D_{t+1}-\hat{D}_t)^2\right]}_{\textup{mutation variance}}.  \label{eq:doob_var}
\end{split}
\end{align}
\end{proposition}
\begin{proof}
It is straightforward to check 
that all the martingale differences $D_{t+1}-\hat{D}_t$ and $\hat{D}_t-D_t$ are uncorrelated with each other and with $D_0$. The proof is finished by writing $D_{T-1} = \sum_{t=0}^{T-1}\sum_{i=1}^N \omega_t^i f(\xi_t^i)$ as a telescoping sum of the martingale differences, computing ${\mathbb E}[D_{T-1}^2]$ in terms of the martingale differences, and subtracting ${\mathbb E}[D_{T-1}]^2 = {\mathbb E}[D_0]^2$ from the resulting expression. 
\end{proof}

In the proof of Theorem~\ref{thm:variance_expressions} below, we will see that the initial condition variance, selection variance, and mutation variance in~\eqref{eq:doob_var} are 
the same as those in~\eqref{eq:var_ic}-\eqref{eq:var_mut}.

\subsection{Proof of the variance formulas}

\begin{proof}[Proof of Theorem~\ref{thm:variance_expressions}]
Using the formula~\eqref{eq:doob_martingale} for the Doob martingale, together with the one-step mean formula~\eqref{eq:one_step1}, the weight update formula~\eqref{eq:weights}, and the conditional independence of particle evolution in~\eqref{eq:mutation},
\begin{align}\begin{split}\label{eq:mut_var_long}
     {\mathbb E}[(D_{t+1} - \hat{D}_{t})^2|\hat{\mathcal F}_{t}] &= \textup{Var}\left.\left(\sum_{i=1}^N \omega_{t+1}^i h_{t+1,T}(\xi_{t+1}^i)\right|\hat{\mathcal F}_t\right) \\
     &= \sum_{i=1}^N (\hat{\omega}_t^i)^2\textup{Var}_{K}h_{t+1,T}(\hat{\xi}_t^i).
     \end{split}
\end{align}
By \eqref{eq:mut_var_long}, the weight update formula~\eqref{eq:hat_weights}, and the bin mean formula~\eqref{eq:bin_mean},
\begin{align*}
    {\mathbb E}[(D_{t+1}-\hat{D}_t)^2] 
    &= {\mathbb E}\left[\sum_{u}  \frac{\omega_t(u)}{N_t(u)} \sum_{i:\hat{\xi}_t^i \in u}\hat{\omega}_t^i \textup{Var}_{K}h_{t+1,T}(\hat{\xi}_t^i)\right] \\
    &= {\mathbb E}\left[\sum_{u}  \frac{\omega_t(u)}{N_t(u)} \sum_{i:{\xi}_t^i \in u}\omega_t^i\textup{Var}_{K}h_{t+1,T}(\xi_t^i)\right].
\end{align*}
In light of Proposition~\ref{lem:doob}, this gives~\eqref{eq:var_mut}. 
Now we turn to~\eqref{eq:var_sel}. 
Using the formula~\eqref{eq:doob_martingale} for the Doob martingale, the one-step mean formula~\eqref{eq:one_step2}, and the fact that selections in distinct bins are conditionally independent,
\begin{align}\begin{split}\label{eq_long_var_sel}
    {\mathbb E}[(\hat{D}_t - D_t)^2|{\mathcal F}_t] 
    &= \textup{Var}\left.\left(\sum_{i=1}^N \hat{\omega}_{t}^i Kh_{t+1,T}(\hat{\xi}_{t}^i)\right|{\mathcal F}_t\right)\\ 
    &= \sum_u \textup{Var}\left.\left(\sum_{i:\hat{\xi}_t^i\in u} \hat{\omega}_t^i Kh_{t+1,T}(\hat{\xi}_t^i)\right|{\mathcal F}_t\right).
    \end{split}
\end{align}
Using~\eqref{eq_long_var_sel}, the weight update formula~\eqref{eq:hat_weights}, and the conditional independence property of multinomial resampling (see {\em e.g.} \cite{douc1}, equation (6)), 
\begin{align*}
    {\mathbb E}[(\hat{D}_t - D_t)^2|{\mathcal F}_t]  
    &= \sum_u \omega_t(u)^2 \textup{Var}\left.\left(\frac{1}{N_t(u)}\sum_{i:\hat{\xi}_t^i\in u} Kh_{t+1,T}(\hat{\xi}_t^i)\right|{\mathcal F}_t\right) \\
    &= \sum_u \frac{\omega_t(u)^2}{N_t(u)}\textup{Var}_{\eta_t^u}(Kh_{t+1,T}).
\end{align*}
Taking expectations in this expression and appealing to Proposition~\ref{lem:doob} gives~\eqref{eq:var_sel}.
\end{proof}

\subsection{Remarks on the variance formulas}

\begin{remark}\label{rmk:total_weight}
We briefly comment on the variance for other methods of the selection and mutation type~\eqref{eq:sel_and_mut}. Consider a method with 
the same mutation step, but a different selection step 
that is still unbiased in the sense of~\eqref{eq:one_step2}. 

In this case, the same variance decomposition~\eqref{eq:doob_var} applies, with selection variance
\begin{equation}\label{eq:gen_sel_var}
\frac{1}{T^2}\sum_{t=0}^{T-2}{\mathbb E}\left[(\hat{D}_t-D_t)^2\right] =\frac{1}{T^2}\sum_{t=0}^{T-2}{\mathbb E}\left[ \textup{Var}\left(\left.\sum_{i=1}^N \hat{\omega}_t^i Kh_{t+1,T}(\hat{\xi}_t^i)\right|{\mathcal F}_t\right)\right].
\end{equation}

For a method like weighted ensemble in which the total weight is always $1$, $Kh_{t+1,T}$ 
can be replaced with $Kh_{t+1,T} - (T-t-1)\int f\,d\mu$ in~\eqref{eq:gen_sel_var} without otherwise
changing the equation. Assumption~\ref{A1} 
shows that $Kh_{t+1,T} - (T-t-1)\int f\,d\mu$ 
is uniformly bounded in $t$ and $T$. As a result, the selection variance~\eqref{eq:gen_sel_var} is $O(1/T)$, and the ergodic theorem 
remains valid.

We observed numerically that if the total weight varies at each time, then the weights tend to approach zero, and the variance~\eqref{eq:WE_var} is of order $T$, as $T \to \infty$. Indeed, $Kh_{t+1,T}$ is typically of order $T$ as $T \to \infty$, which suggests that in this case, the selection variance~\eqref{eq:gen_sel_var} is on the order of $\sum_{t=0}^{T-2}{\mathbb E}[\textup{Var}(\sum_{i=1}^N \hat{w}_t^i|{\mathcal F}_t)]$ as $T \to \infty$. Of course, the ergodic theorem fails if the variance~\eqref{eq:WE_var} goes to infinity as $T \to \infty$.
\end{remark}

\begin{remark}\label{rmk:DMC}
Note that direct Monte Carlo -- which we define as independent, equally-weighted particles -- is a special case of weighted ensemble in which each particle always has weight $1/N$, every parent always gets its own bin $u$, and $N_t(u)$ is always $1$. In this case the selection variance is zero, 
while the mutation variance is
\begin{equation}\label{eq:DMC_mutvar}
\frac{1}{NT^2} \sum_{t=0}^{T-2}{\mathbb E}\left[\frac{1}{N}\sum_{i=1}^N \textup{Var}_Kh_{t+1,T}(\xi_t^i)\right] = \frac{1}{NT^2} \sum_{t=0}^{T-2}\nu(K^t(\textup{Var}_Kh_{t+1,T})).
\end{equation}
\end{remark}

\begin{remark}\label{rmk:residual}
If the selections in the bins use residual multinomial resampling~\cite{douc1} instead of multinomial resampling, then the selection variance~\eqref{eq:var_sel} becomes
\begin{equation*}
\frac{1}{T^2}\sum_{t=0}^{T-2} {\mathbb E}\left[\sum_u\frac{\omega_t(u)^2}{N_t(u)^2}r_t(u)\textup{Var}_{\gamma_t^u} (Kh_{t+1,T})\right],
\end{equation*}
where $r_t(u)$ and $\gamma_t^u$ are defined from the residuals $r_t^i = \frac{N_t(u)\omega_t^i}{\omega_t(u)} - \big\lfloor \frac{N_t(u)\omega_t^i}{\omega_t(u)}\big\rfloor$ by $$r_t(u) = \sum_{i:\xi_t^i \in u} r_t^i, \qquad \gamma_t^u = \sum_{i:\xi_t^i\in u} \frac{r_t^i}{r_t(u)}\delta_{\xi_t^i}.$$
We omit proof, but include this formula in case 
it is useful for optimizations, since residual multinomial resampling performs much better than 
multinomial resampling, yet still admits simple explicit variance expressions.
\end{remark}

\subsection{Proof of the ergodic theorem}

\begin{lemma}\label{lem:var_scaling}
If Assumption~\ref{A1} holds, then as $T \to \infty$, 
\begin{equation*}
\textup{Var}\left(\frac{1}{T}\sum_{t=0}^{T-1}\sum_{i=1}^N \omega_t^i f(\xi_t^i)\right) = O(1/T).
\end{equation*}
\end{lemma}

\begin{proof}
By Assumption~\ref{A1}, we have 
\begin{equation}\label{geom_erg}
|K^tf(x)-K^tf(y)| \le c\lambda^t, \qquad \text{for all }x,y\text{ and all }t\ge 0,
\end{equation} 
where now $c>0$ is a different constant. Thus
\begin{equation*}
|h_{t+1,T}(x)-h_{t+1,T}(y)| \le \sum_{s=0}^{T-t-2}|K^sf(x)-K^sf(y)| \le \frac{c}{1-\lambda} := C.
\end{equation*}
This shows that $\textup{Var}_\eta h_{t+1,T} \le C^2$, and similarly
$\textup{Var}_\eta Kh_{t+1,T} \le C^2$, for any 
probability distribution $\eta$. As a result, 
the selection and mutation variances~\eqref{eq:var_sel}-\eqref{eq:var_mut} are both $O(1/T)$ as $T \to \infty$. 
By similar arguments,  
$\textup{Var}(\sum_{i=1}^N \omega_0^i h_{0,T}(\xi_0^i)) \le C^2$, which makes
the initialization variance~\eqref{eq:var_ic} $O(1/T^2)$. 
Thus, the variance in~\eqref{eq:WE_var} is $O(1/T)$.
\end{proof}

\begin{proof}[Proof of Theorem~\ref{thm_main}]
Define $\theta_T = \frac{1}{T}\sum_{t=0}^{T-1}\sum_{i=1}^N \omega_t^i f(\xi_t^i)$. Then for $0 < S \le T$,
\begin{equation}\label{theta_diff}
|\theta_T - \theta_S| = \left|\left(\frac{S}{T}-1\right)\theta_S + \frac{1}{T} \sum_{t=S}^{T-1}\sum_{i=1}^N\omega_t^i f(\xi_t^i)\right| \\ 
\le 2\sup|f|\left(1-\frac{S}{T}\right).
\end{equation}
By Lemma~\ref{lem:var_scaling}, $\textup{Var}(\theta_T) = O(1/T)$. So by Chebyshev's inequality, there is $c>0$ so that 
\begin{equation*}
{\mathbb P}\left(\left|\theta_T - {\mathbb E}[\theta_T]\right| \ge cT^{1/3-1/2}\right) \le \frac{1}{T^{2/3}}
\end{equation*}
for large enough $T$. 
With $T_n = n^2$, by the Borel-Cantelli lemma, 
there is $n_0$ such that 
\begin{equation}\label{BorelCantelli}
\left|\theta_{T_n} - {\mathbb E}[\theta_{T_n}]\right| < cT_n^{1/3-1/2} \qquad\text{a.s. for all } n \ge n_0.
\end{equation}
By the unbiased property and 
Assumption~\ref{A1}, 
\begin{equation}\label{eq:erg_limit}\lim_{T \to \infty} {\mathbb E}[\theta_T] = \int f\,d\mu.
\end{equation} 
Now given $S>0$, we can choose $T_n$ so that $T_n \le S \le T_{n+1}$ and write
\begin{equation}\label{triangle_ineq}
\left|\theta_S - \int f\,d\mu\right| \le |\theta_S - \theta_{T_n}| + |\theta_{T_n} - {\mathbb E}[\theta_{T_n}]| 
+ \left|{\mathbb E}[\theta_{T_n}] - \int f\,d\mu\right|.
\end{equation}
By~\eqref{theta_diff}-\eqref{eq:erg_limit}, with probability $1$, the right hand side of~\eqref{triangle_ineq} vanishes as $S \to \infty$.
\end{proof}

\acks 
\noindent The author would like 
to acknowledge Fr{\'e}d{\'e}ric C{\'e}rou, Peter Christman, Josselin Garnier, Gideon Simpson, Gabriel Stoltz and
Brian Van Koten for helpful comments, 
and especially Jeremy Copperman, Matthias Rousset, Robert J. Webber,
and Dan Zuckerman for interesting discussions 
and insights. 
The author thanks 
Robert J. Webber for pointing out 
errors and making many helpful suggestions 
concerning a previous version of the 
manuscript.

\fund 
\noindent The author also gratefully acknowledges support from the National Science Foundation via the awards NSF-DMS-1818726 and
NSF-DMS-1522398.

\end{document}